\documentclass[titlepage,11pt]{article}
\usepackage{latexsym}
\usepackage{amsfonts}
\usepackage{amsmath}
\usepackage[T1]{fontenc}
\newcommand\blackslug{\hbox{\hskip 1pt \vrule width 4pt height 8pt depth 1.5pt
        \hskip 1pt}}
\newcommand\bbox{\hfill \quad \blackslug \bigbreak}
\newcommand{\vare}{\varepsilon}
\def\DD{\hbox{-}}
\def\CC{\hbox{-}\cdots\hbox{-}}
\def\LL{,\ldots,}
\def\cupcup{\cup\cdots\cup}

\title{Polynomial bounds for chromatic number \\ VIII. Excluding a path and a complete multipartite graph}
\author{Tung Nguyen\thanks{Supported by AFOSR grant FA9550-22-1-0234 and NSF
grant DMS-2154169.}\\
Princeton University, Princeton, NJ 08544
\\
\\
Alex Scott\thanks{Research supported by EPSRC grant EP/X013642/1.}\\
Mathematical Institute, University of Oxford, Oxford OX2 6GG, UK
\\
\\
Paul Seymour\thanks{Supported by AFOSR grant FA9550-22-1-0234 and NSF
grant DMS-2154169.}\\
Princeton University, Princeton, NJ 08544}

\date{January 9, 2023; revised \today}

\newtheorem{thm}{}[section]

\newcommand{\Proof}{\noindent{\bf Proof.}\ \ }

\begin{document}
\maketitle
\begin{abstract}

We prove that for every path $H$, and every integer $d$, there is a polynomial $f$ such that every graph $G$ with chromatic number greater than $f(t)$
either contains $H$ as an induced subgraph, or contains as a subgraph the complete $d$-partite graph with parts of 
cardinality $t$.  For $t=1$ and general $d$ this is a classical theorem of 
Gy\'arf\'as, and for $d=2$ and
general $t$ this is a theorem of Bonamy et al. 


\end{abstract}

\section{Introduction}

A graph is {\em $H$-free} if it contains no induced copy
of $H$.
The Gy\'arf\'as-Sumner conjecture~\cite{gyarfastree,sumner} says:
\begin{thm}\label{GSconj}
{\bf Conjecture: }For every forest $H$, and every $k$, every $H$-free graph that does not contain
a clique on $k$ vertices has bounded chromatic number.
\end{thm}

This is known only for a few special kinds of forest (see \cite{survey} for a survey).  It is known when $H$ is a path:
this was proved by Gy\'arf\'as~\cite{gyarfastree} in the original paper proposing the Gy\'arf\'as-Sumner conjecture, with an
elegant proof technique now known as the ``Gy\'arf\'as path'' method. But even when $H$ is a path, and indeed even when $H$
is the five-vertex path $P_5$, it is not known whether the bound in \ref{GSconj} is polynomial in $k$ (the best current bound is 
that every $P_5$-free graph $G$ which does not contain a clique on $k$ vertices has chromatic number 
$\chi(G)\le k^{\log_2k}$~\cite{poly4}).

There is a stronger conjecture, which has a complicated parentage (it resulted from the joining of \ref{GSconj} and an open 
special case of a false conjecture of Esperet~\cite{esperet}):
\begin{thm}\label{esperet}
{\bf Conjecture: }For every forest $H$, there is a polynomial $f(k)$ such that for every $k\ge 0$,
every $H$-free graph that does not contain
a clique on $k$ vertices has chromatic number at most $f(k)$.
\end{thm}
So far the conjecture is only known for a few forests $H$, including all trees that do not contain $P_5$ as a subgraph (see \cite{poly2,poly3,poly6,rs}).  

If we exclude a complete bipartite graph $K_{t,t}$ rather than a clique, we can do better: 
the conjectures \ref{GSconj} and \ref{esperet} become theorems. R\"odl~\cite{rodl} proved the analogue of \ref{GSconj}, 
that:
\begin{thm}\label{rodl}
For every forest $H$, and every $t$, every $H$-free graph that does not contain
$K_{t,t}$ as a subgraph has bounded chromatic number.
\end{thm}
The analogue of \ref{esperet} (that is, \ref{rodl} with polynomial bounds) is also now known.  It was proved first for the case when $H$ is a path by
Bonamy et al.~\cite{bonamy}. 
Two of us, with Spirkl~\cite{poly1}, subsequently
proved it for all forests:
\begin{thm}\label{bonamy}
For every forest $H$ there is a polynomial $f(t)$ such that for every $t\ge 1$, every $H$-free graph that does 
not contain $K_{t,t}$ as a subgraph
has chromatic number at most $f(t)$.
\end{thm}

What about excluding the {\em tripartite} graph $K_{t,t,t}$? Or, more generally, excluding as a subgraph the complete $d$-partite graph with parts of cardinality $t$, which we 
denote by $K_d(t)$? 
We will show:
\begin{thm}\label{paththm}
For every path $H$ and all $d\ge 1$, there is a polynomial $f$ such that,  for all $t\ge 1$, if a graph $G$ is $H$-free and does not
contain $K_d(t)$ as a subgraph, then $\chi(G)\le f(t)$.
\end{thm}
For $t=1$ and general $d$ this is Gy\'arf\'as' result for paths, and for $d=2$ and
general $t$ this is the theorem of Bonamy et al. 

Let us say a graph $H$ is {\em multibounding} if
for every $d\ge 1$, there is a polynomial $f(t)$ such that, for all $t\ge 1$,
if a graph $G$ is $H$-free and does not contain $K_d(t)$ as a subgraph, then $\chi(G)\le f(t)$.
It is easy to see that every multibounding graph $H$ must be a forest and must 
satisfy the Gy\'arf\'as-Sumner conjecture (take $t=1$), and any forest satisfying \ref{esperet} is multibounding.  We propose:
\begin{thm}\label{ourconj}
{\bf Conjecture: } Every forest is multibounding.
\end{thm}
Among all the forests that are known to satisfy \ref{GSconj}, which can we show are multibounding? 
\begin{itemize}
\item Two of us~\cite{poly5} showed that 
every tree of radius two is multibounding, strengthening the theorem of Kierstead and Penrice~\cite{kierstead} 
that every tree $H$ of radius two satisfies \ref{GSconj}.
\item In this paper we will show that paths are multibounding, strengthening Gy\'arf\'as' theorem that paths satisfy \ref{GSconj}.
\item A {\em broom} is obtained from a path with one end $v$ by adding leaves adjacent to $v$, and we will also show 
that brooms are multibounding, strengthening a result of Gy\'arf\'as that brooms satisfy \ref{GSconj}.
\item Finally, we will show that if $H$ is a forest and each of its components is multibounding
then $H$ is multibounding.
\end{itemize}

\section{The proof of \ref{paththm}}

The proof uses a new version of the ``Gy\'arf\'as path'' method.  The usual Gy\'arf\'as path technique goes by induction on the 
clique number: the idea (very roughly) is that the inductive hypothesis allows us to assume that neighbourhoods have small 
chromatic number, and it is then possible to walk step-by-step into the graph, always heading towards an unexplored portion of 
the graph with large chromatic number.  At each step, we have visited only a bounded number of vertices, and so the portion of 
the graph in the neighbourhood of the explored set has bounded chromatic number.

This argument does not give polynomial bounds in our setting, as the induction pushes up the bounds  too quickly.  Here, we are 
concerned with complete $d$-partite graphs $K_d(t)$, and argue by induction on $d$.  Rather than a simple induction using 
neighbourhoods, we instead use the fact that for any $t$ vertices their {\em common} neighbourhood does not contain $K_{d-1}(t)$ 
and so by induction has small chromatic number (at most polynomial in $t$).  It follows that, if $H$ is a subgraph with large 
chromatic number, then there are at most $t-1$ vertices adjacent to almost all (in an appropriate sense) of $H$.  We now make an 
argument analogous to the Gy\'arf\'as path argument, but rather than one step at a time, we move in a sequence of short dashes, 
and need to be more careful about the parts of the graph that we explore. At each stage, we make sure we are attached to a 
subgraph $H$ that is both highly connected and has large chromatic number; we identify the next subgraph $H'\subseteq H$ (again, highly 
connected and with large chromatic number), and then take a path through $H$ to $H'$, using the connectivity to 
avoid any (at most $t-1$) ``dangerous'' vertices that can see almost all of $H'$.  Iterating the argument gives the required path.

We need a lemma:
\begin{thm}\label{conn}
For all integers $t>0$, if $G$ has chromatic number at least $4t$, then 
$G$ has a $t$-connected induced subgraph with chromatic number at least $\chi(G)-2t$.
\end{thm}
This is proved in~\cite{tung}, but weaker statements have been proved in several other papers: for instance, by Gir\~ao and 
Narayanan~\cite{narayanan}, and by Penev, Thomass\'e and Trotignon~\cite{thomasse}, who proved:
\begin{thm}\label{stephan}
For all integers $t>0$, if $G$ has chromatic number at least $2(t-1)^2+1$, then
$G$ has a $t$-connected induced subgraph with chromatic number at least $\chi(G)-2t+3$.
\end{thm}
(Any statement 
of this form is good enough for the application here, with appropriate adjustment of the constants.)

If $X\subseteq V(G)$, we define $\chi(X)=\chi(G[X])$.
Let $X,Y$ be disjoint subsets of $V(G)$, where $|X|=t$ and every vertex in $Y$ is adjacent to every vertex in $X$. We call
$(X,Y)$ a {\em $t$-biclique}, and its {\em value} is $\chi(Y)$.
For integers $p,t\ge 1$, let us say a {\em $(p,t)$-balloon $(P,Y)$} in $G$ consists of an induced path $P$ of $G$ with vertices
$v_1\DD v_2\CC v_p$ in order, and a subset $Y\subseteq V(G)$,
 with the following properties:
\begin{itemize}
\item $v_1\LL v_{p-1}\notin Y$, $v_p\in Y$,
none of $v_1\LL v_{p-2}$ have a neighbour in $Y$, and if $p\ge 2$ then
$v_p$ is the unique neighbour of $v_{p-1}$ in $Y$; and
\item $G[Y]$ is $t$-connected.
\end{itemize}
The {\em value} of the $(p,t)$-balloon $(P,Y)$ is $\chi(Z)$, where $Z$ is the set of vertices in $Y$ nonadjacent to $v_p$.

\begin{thm}\label{bigt}
For every graph $G$, and all integers $p,q,s,t\ge 1$, if $G$ contains no  $(p,t)$-balloon of value at least $q$, and
$G$ contains no $t$-biclique of value $s$, then
$$\chi(G)\le (1+t+\cdots+t^{p-1})(s+t(2t+9)) + t^pq.$$
\end{thm}
\Proof
For integers $p,q,s,t\ge 1$, define
$$k(p,q,s,t) = (1+t+\cdots+t^{p-1})(s+t(2t+9)) + t^pq.$$
We proceed by induction on $p$. Suppose first that $p=1$, and $\chi(G)>k(1,q,s,t)$. 
We want to find a $(1,t)$-balloon of value
at least $q$; that is, a set $Y$ and a vertex $v\in Y$ such that $G[Y]$ is
$t$-connected and the non-neighbours of $v$ in $Y$ induce a subgraph with chromatic
number at least $q$.
By \ref{conn}, since $\chi(G)>k(1,q,s,t)\ge 4t$,
$G$ has a $t$-connected induced subgraph  $G[Z]$ say with chromatic number at least $\chi(G)-2t$.
Since $\chi(G)-2t\ge k(1,q,s,t)-2t\ge t$,
there exist $t$ distinct vertices in $Z$, say $z_1\LL z_t$. For each $i$, let 
$X_i$ be the set of vertices in $Z\setminus \{z_i\}$ nonadjacent to $z_i$.
Since $G[Z]$
is $t$-connected,  we may assume 
that $G[X_i]$ has chromatic number at most $q$ for $1\le i\le t$, and hence so does $G[X_i\cup \{z_i\}]$ since $q\ge 1$. 
Thus the union
of the sets $X_1\LL X_t$ and $\{z_1\LL z_t\}$ has chromatic number at most $tq$. Moreover, since
$G$ contains no $t$-biclique of value $s$, the set of vertices
in $Z\setminus \{z_1\LL z_t\}$ adjacent to all of $z_1\LL z_t$ has chromatic number less than $s$.
Consequently $\chi(Z)\le s+tq$, and so $\chi(G)\le s+t(q+2)\le k(1,q,s,t)$ and the theorem holds.

Thus we may assume that $p\ge 2$, and theorem holds for $p-1$ (and all values of $q,s,t$).  Let $q'=s+t(2t+9)+tq$.
In particular, we may assume that
$G$ contains a $(p-1,t)$-balloon of value $q'$, since $k(p-1,q',s,t)= k(p,q,s,t)$.
We assume for a contradiction that $\chi(G)>k(p,q,s,t)$.
Let $(P',Y')$ be a $(p-1,t)$-balloon of value $q'$, where $P'$ has vertices 
$v_1\CC v_{p-1}$. 
Let $Z'$ be the set of vertices in $Y'$ nonadjacent to $v_{p-1}$; thus, $\chi(Z')\ge q'$. 
Since $q'\ge 12t$,
\ref{conn} implies that there exists $Y\subseteq Z'$ such that $G[Y]$ is $3t$-connected and $\chi(G[Y])\ge q'-6t$.
Let $X$ be the set of all vertices $v\in Y'$ such that the set of vertices in $Y$ different from and nonadjacent to $v$
has chromatic number less than $q+2t$. 
\\
\\
(1) {\em $|X|<t$.}
\\
\\
Suppose that $|X|\ge t$, and choose $X'\subseteq X$ with $|X'|=t$. Since $G$ contains no $t$-biclique of value $s$, the set
of vertices in $Y\setminus X'$ that are adjacent to every vertex in $X'$ has chromatic number less than $s$. 
But for each $v\in X'$, the set of vertices in $Y$ different from and nonadjacent to $v$
has chromatic number less than $q+2t$, and hence so does the union of this set and $\{v\}$. But $Y$ is the union of
these subsets; so $\chi(Y)<s+t(q+2t)\le q'-6t\le \chi(Y)$, a contradiction. This proves (1).

\bigskip
Since $G[Y']$ is $t$-connected, and $v_{p-1}\notin X$, and $Y\not\subseteq X$ (because $G[Y]$ is $3t$-connected and hence 
$|Y|\ge 3t+1\ge t$), there is an induced path $Q$ of $G[Y']$ between $v_{p-1}$ and $Y$ with no vertex in $X$. 
Choose $Q$ minimal, and let its vertices be $v_{p-1}=u_1\CC u_k$ in order, where $u_k\in Y$. Hence $k\ge 3$, since
$u_1=v_{p-1}$ has no neighbour in $Y$. From the minimality of $Q$, 
none of $u_1\LL u_{k-1}$
belong to $Y$, and none of $u_1\LL u_{k-2}$ have a neighbour in $Y\setminus X$. Let $R$ be the path with vertices
$v_1\CC v_{p-1} \DD u_2\CC u_{k-1}\DD u_k$. Now there are two cases.

First, suppose that $u_{k-1}$ has at least $t$ neighbours in $Y\setminus X$. Let $R'$ be a subpath of $R\setminus \{u_k\}$ with $p$ vertices 
and one end $u_{k-1}$. Since $G[(Y\setminus X)\cup \{u_{k-1}\}]$ is $t$-connected (because $G[Y]$ is $3t$-connected, and $|X|<t$, 
and $u_{k-1}$ has at least $t$ neighbours in $Y\setminus X$),
and the set of non-neighbours of $u_{k-1}$ in $Y\setminus X$ has chromatic number at least $q+t\ge q$ (because 
$u_{k-1}\notin X$ and $|X|<t$) it follows that $(R', (Y\setminus X)\cup \{u_{k-1}\})$ is a $(p,t)$-balloon of value at least $q$,
a contradiction. 

Now suppose that $u_{k-1}$ has fewer than $t$ neighbours in $Y\setminus X$. The vertex $u_k$ is one such neighbour;
let $N$ be the set of the neighbours of $u_{k-1}$ in $Y\setminus X$ different from $u_k$. Let $R'$ be a subpath of $R$
with $p$ vertices and one end $u_k$. Since $G[Y\setminus (X\cup N)]$ is $t$-connected (because $G[Y]$ is $3t$-connected and 
$|X\cup N|\le 2t$),
and the set of non-neighbours of $u_k$ in $Y\setminus (X\cup N)$ has chromatic number at least $q$ (because 
$u_k\notin X$ and $|X\cup N|<2t$) it follows that $(R', Y\setminus (X\cup N))$ is a $(p,t)$-balloon of value at least $q$,
a contradiction. This proves \ref{bigt}.~\bbox

Let us deduce that all paths are multibounding, that is, \ref{paththm}, which we restate:
\begin{thm}\label{paththm2}
For all $d\ge 1$ and every path $H$, there is a polynomial $f$ such that,  for all $t\ge 1$, if a graph $G$ is $H$-free and does not
contain $K_d(t)$ as a subgraph, then $\chi(G)\le f(t)$.
\end{thm}
\Proof
Let $H$ have $p$ vertices. $H$-free graphs contain no $(p,t)$-balloon of value $1$, so by \ref{bigt}, 
there exist integers $b,c\ge 1$ such that for all integers $s,t\ge 1$, and every graph $G$,
if $G$ is $H$-free and contains no $t$-biclique of value $s$, then $\chi(G)\le bt^cs$.

We prove by induction on $d$ that if $G$ is $H$-free and does not
contain $K_d(t)$ as a subgraph, then $\chi(G)\le (bt^c)^d$. If $d=1$, graphs not containing $K_d(t)$ as a subgraph have fewer
than $t$ vertices, and so
have chromatic number less than $t\le bt^c$ (since $b,c\ge 1$) and the claim holds.
Thus we assume that $d>1$ and the claim holds
for $d-1$.
Let $G$ be $H$-free, with chromatic number more than $(bt^c)^d$; we will show that $G$
contains $K_d(t)$ as a subgraph.
Since $\chi(G)>(bt^c)^d$, we deduce from the choice of $b,c$ (taking $s=(bc^t)^{d-1}$)
that $G$
contains a $t$-biclique $(X,Y)$ of value $(bc^t)^{d-1}$.
From the inductive hypothesis applied to $G[Y]$,
$G[Y]$ contains $K_{d-1}(t)$ as a subgraph. But then, adding $X$, we find that
$G$ contains $K_d(t)$ as a subgraph. This proves \ref{paththm}.~\bbox

The proof we just gave is redundant, because \ref{paththm} follows from the analogous theorem \ref{broomthm}
for brooms which we will prove later,
but the proof above is a little simpler than the proof of \ref{broomthm}.

\section{Brooms}

A {\em broom} is obtained from a path with one end $v$ by adding leaves adjacent to $v$. In this section we deduce from
\ref{bigt} that brooms are multibounding.
Essentially the proof of \ref{paththm} still works.
For $k\ge 1$, let $S_k$ be the star with $k+1$ vertices: that is, a tree in which some vertex is incident with every edge.
Let $\omega(G)$ denote the size of the largest clique of $G$.
We need the following lemma, due to Gy\'arf\'as~\cite{gyarfas}:
\begin{thm}\label{nostar}
If $G$ is $S_k$-free, then $\chi(G)\le \omega(G)^k$.
\end{thm}
We deduce that brooms are multibounding, that is:

\begin{thm}\label{broomthm}
For all $d\ge 1$ and every broom $H$, there is a polynomial $f$ such that,  for all $t\ge 1$, if a graph $G$ is $H$-free and does not
contain $K_d(t)$ as a subgraph, then $\chi(G)\le f(t)$.
\end{thm}
\Proof
Let $H$ be obtained from a $p$-vertex path by adding $r$ new vertices each adjacent to one (and the same) end of the path.
We prove the result by induction on $d$.

If $d=1$ then the result holds trivially, because graphs without $K_1(t)$ as a subgraph have fewer than $t$
vertices.
So we assume that $d>1$ and the result holds for $d-1$. Let $g(t)$ be the corresponding polynomial.
Define
$$f(t) = (1+t+\cdots+t^{p-1})(g(t)+1+t(2t+9)) + t^p(dt)^{2r}.$$
Now let $G$ be an $H$-free graph that does not contain $K_d(t)$ as a subgraph; we will show that $\chi(G)\le f(t)$.
\\
\\
(1) {\em $G$ contains no $(p,t)$-balloon of value $(dt)^{2r}$.}
\\
\\
Suppose that $(P,Y)$ is such a $(p,t)$-balloon, where $P$ has vertices $v_1\CC v_p$ in order.
Since $G$ does not contain $K_d(t)$ as a subgraph, it follows that $\omega(G)< dt$, and so $\chi(Y)> \omega(G)^{2r}$,
and therefore $G[Y]$ contains an induced copy of $S_{2r}$ by \ref{nostar}.
Let $a\in Y$, and let $\{b_1\LL b_{2r}\}\subseteq Y\setminus \{a\}$ be a stable set of neighbours of $a$.
Since $G[Y]$ is connected, there is an induced path $Q$ from $v_p$  to $a$, say with vertices $v_p=q_1\CC q_n=a$.
Choose $i\in \{1\LL n\}$  minimum such that $q_i$ belongs to or is adjacent to a vertex in $\{a,b_1\LL b_{2r}\}$.
If $q_i\in \{a,b_1\LL b_{2r}\}$, then $i = 1$ from the minimality of $i$, and the subgraph induced on
$\{v_1\LL v_p,a,b_1\LL b_{2r}\}$ contains a copy of $H$, a contradiction. So $q_i\ne \{a,b_1\LL b_{2r}\}$. Let $N$
be the set of neighbours in $\{a,b_1\LL b_{2r}\}$ of $q_i$. If $N$ contains at least $r$ of $\{b_1\LL b_{2r}\}$, then the
subgraph induced
on $\{v_1\LL v_p,q_1\LL q_i\}\cup N$ contains a copy of $H$, a contradiction. So we may assume that $b_1\LL b_r\notin N$.
If $a\in N$ then the subgraph induced
on $\{v_1\LL v_p,q_1\LL q_i, a, b_1\LL b_r\}$
contains a copy of $H$, a contradiction. If $a\notin N$, we may assume that $b_{r+1}\in N$, and then the
subgraph induced
on $\{v_1\LL v_p,q_1\LL q_i, a, b_1\LL b_r, b_{r+1}\}$
contains a copy of $H$, a contradiction. This proves (1).
\\
\\
(2) {\em $G$ contains no $t$-biclique of value $g(t)+1$.}
\\
\\
If $(X,Y)$ is such a $t$-biclique, from the inductive hypothesis $G[Y]$ contains $K_{d-1}(t)$ as a subgraph, and hence
$G$ contains $K_d(t)$ as a subgraph, a contradiction. This proves (2).

\bigskip

From (1), (2) and \ref{bigt},
$$\chi(G)\le (1+t+\cdots+t^{p-1})(g(t)+1+t(2t+9)) + t^p(dt)^{2r}=f(t).$$
This proves \ref{broomthm}.~\bbox

For every forest $H$ that is known to satisfy the Gy\'arf\'as-Sumner conjecture \ref{GSconj}, one might hope to show that $H$
is multibounding. So far, we know that trees of radius two, and brooms, are multibounding. 
What is the next simplest forest we could try? There are several candidates for this: disjoint unions of trees that are 
known to be multibounding? Subdivisions of a star, or of a claw? Double brooms? 

We have answered the first of these, as we shall explain in the next section, and 
made some progress on the third, as we explain now. A {\em double broom}
is a tree obtained from a path by adding new vertices each adjacent to an end of the path, but not necessarily all to the same end.
It was proved in~\cite{brooms}  that double brooms satisfy \ref{GSconj}.

Let $B(p,r,s)$ be the double broom obtained from a $p$-vertex path by adding $r$ leaves at one end and $s$ at the other.
We have not been able to show that double brooms are multibounding, but we proved the following:
\begin{thm}\label{doublebroom}
For all $d,p,r,s,t\ge 1$, there is a polynomial $f(t)$, such that if a graph $G$ does not contain $K_d(t)$, and its chromatic number 
is more than $f(t)$, then $G$ contains a double broom $B(p',r,s)$ for some $p'\ge p$.
\end{thm}

The proof is a modification of the proof of \ref{bigt}, beginning with a large star induced subgraph and growing 
a balloon tethered to this star.
Here is a sketch. We may assume that $r=s$, and that $G$ has chromatic number at least some large polynomial in $t$, and $G$ does not 
contain a $t$-biclique with large value (again, value at least a large polynomial in $t$, but not so large). We may assume that $G$
is $t$-connected, by \ref{conn}. For each $v\in V(G)$, $N(v)$ and $M(v)$ denote its set of neighbours and non-neighbours.
Since there is no $t$-biclique with large value, there are fewer than $t$ vertices $v$ for which $\chi(M(v))$ is small. So there
exists $v$ with $|N(v)|$ and $\chi(M(v))$ large. Choose $N\subseteq N(v)$ with cardinality equal to a large polynomial in $t$.

Suppose there is a stable subset
$S\subseteq N$ with $|S|=3s$ such that the set $X$ of vertices in $M(v)$ with no neighbour in $S$ has large chromatic number. Apply
\ref{conn} to obtain a $t$-connected subset $Y$ of $X$, still with large $\chi$. Fewer than $t$
vertices in $V(G)\setminus Y$ are adjacent to almost all of $Y$, so there is a path from $v$ to $Y$ avoiding such vertices.
Let $v,v_1,v_2$ be the first three vertices of this path (we may assume that no other vertices of the path have neighbours in $S$).
Every vertex in $S$ is either adjacent to neither of $v_1,v_2$, or adjacent to $v_1$ and not to $v_2$, or adjacent to $v_2$;
and since $|S|=3s$, at least $s$ vertices in $S$ behave the same. In each case, 
we have half the double broom, and now use the argument of \ref{broomthm} to obtain the other
half. 

If there is no such $S$, then since there are only polynomially many choices of stable sets in $N$ with size $3s$ (because $|N|$
is a polynomial in $t$),
the set of all vertices 
$u\in M(v)$ such that $N(v)\cap M(u)$ contains a stable set of size $3s$ has small chromatic number; and so the set $Z$ of $u\in M(v)$
such that $N(v)\cap M(u)$ contains no such stable set has large chromatic number. But if $N(v)\cap M(u)$ contains no stable set
of size $3s$, then its cardinality is at most a polynomial in $t$ (by Ramsey's theorem, since $\omega(G)\le td$). So every vertex
in $Z$ is adjacent to almost all of $N(v)$. Inductively, $G[Z]$ contains $K_{d-1}(t)$ as a subgraph ($K$ say), and each vertex of $K$
has only a few non-neighbours in $N(v)$; so $K$ has at least $t$ common neighbours in $N(t)$, giving a copy of $K_d(t)$, 
a contradiction. This ends the sketch of the proof of \ref{doublebroom}.

\section{Disjoint unions}

In this section we show that:
\begin{thm}\label{comps}
If $H$ is a forest and every component of $H$ is multibounding, then $H$ is multibounding.
\end{thm} 

We will deduce this by applying a strengthened version of a theorem from~\cite{pairs}. 
It turns out that a modification of the proof
in~\cite{pairs} gives a stronger result than was noticed at the time. We say $J$ is the {\em complete join} of $J_1,J_2$ if the complement 
of $J$
is the disjoint union of the complements of $J_1$ and $J_2$.

The theorem from~\cite{pairs} says that if $H$ is the disjoint union of $H_1,H_2$, and $J$ is the complete join 
$J_1$ and $J_2$, and $G$ is both $H$- and $J$-free, then $V(G)$ can be partitioned into a bounded number of parts
such that for each part $X$, $G[X]$ is either $H_1$-free or $H_2$-free or $J_1$-free or $J_2$-free (unless all of $H_1,H_2,J_1,J_2$
have only one vertex, when the theorem is false.) We need to improve this in two ways. First, 
we want to tighten the bound on the nunmber of parts, which turns out only to depend on $|H_1|$ and $|J_1|$, and is 
polynomial in $|J_1|$; and second, we want to apply it to \ref{comps}, taking $J$ to be $K_d(t)$. This presents a difficulty,
because in \ref{comps}, we are working with containing $K_d(t)$ as a subgraph, not as an induced subgraph, and we need to work around
this somehow. The easiest fix seems to be to work with excluding sets of graphs rather than single graphs: because excluding, say,
$K_{t,t}$
as a subgraph is equivalent to excluding as induced subgraphs all graphs with $2t$ vertices in which half the vertices are adjacent
to the other half.

Let $\mathcal{H}$ be a nonempty class of graphs. A graph is {\em $\mathcal{H}$-free} if it is $H$-free for each 
$H\in \mathcal{H}$. The {\em max-size} of $\mathcal{H}$ is the maximum number of vertices of members of $\mathcal{H}$, if this maximum 
exists, and $\infty$ otherwise; and the {\em min-size} is the minimum number of vertices of members of $\mathcal{H}$.
If $\mathcal{H}_1, \mathcal{H}_2$ are classes of graphs, we denote by $\mathcal{H}_1\uplus \mathcal{H}_2$
the class of all graphs that are the disjoint union of a member of $\mathcal{H}_1$ and a member of $\mathcal{H}_2$; and
$\mathcal{H}_1 * \mathcal{H}_2$ is the class of all graphs that are the complete join of a member of $\mathcal{H}_1$
and a member of $\mathcal{H}_2$.

We will show:
\begin{thm}\label{pairsthm}
Let $\mathcal{H}_1, \mathcal{H}_2, \mathcal{J}_1, \mathcal{J}_2$ be nonempty classes of non-null graphs, and let $k\ge 0$ be 
an integer such that $\mathcal{H}_1$ 
and $\mathcal{J}_1$
both have max-size at most $k-2$.
Let at least one of $\mathcal{H}_1, \mathcal{H}_2, \mathcal{J}_1, \mathcal{J}_2$ have min-size at least two.
Let $S\in \mathcal{H}_1$, let  $s=|S|$, and let $n=sk^s+3$.  If $G$ is $\mathcal{H}_1\uplus\mathcal{H}_2$-free and
$\mathcal{J}_1 * \mathcal{J}_2$-free, then 
$V(G)$ can be partitioned into $n$ sets $X_1\LL X_n$, such that for $1\le i\le n$, $G[X_i]$ is either
$\mathcal{H}_1$-free or $\mathcal{H}_2$-free or  $\mathcal{J}_1$-free or  $\mathcal{J}_2$-free.
\end{thm}
\Proof Let us say a subset $X\subseteq V(G)$ is {\em good} if $G[X]$ is either
$\mathcal{H}_1$-free or $\mathcal{H}_2$-free or  $\mathcal{J}_1$-free or  $\mathcal{J}_2$-free.
For $X\subseteq V(G)$, let $\mu(X)$ be the minimum $t$ such that
there is a partition of $X$ into $t$ good sets. (This exists, since at least one of the four classes has min-size at least two, and
so all one-vertex sets are good.) We call $\mu(X)$ the {\em measure} of $X$.
Let $V(S) = \{h_1\LL h_s\}$, and let us add a new vertex $h_0$ to $S$, with no neighbours, forming a graph $S^+$.

Now let $1\le q\le s$, and assume that $\phi$ is an isomorphism from $S[\{h_q\LL h_s\}]$ to some induced subgraph $F$ of $G$. 
We call $\phi$ a {\em $q$-map}, and 
we say that $\phi$ is {\em $t$-general} if for $0\le p<q$ there is a subset $Y_p$ of $V(G)\setminus V(F)$ with the following properties:
\begin{itemize}
\item the sets $Y_{0}\LL Y_{q-1}$ are pairwise disjoint, and all have measure at least $t$;
\item for $q\le r\le s$, and $0\le p<q$, if $h_r, h_p$ are adjacent in $S^+$ then $\phi(h_r)$ is adjacent to every vertex
of $Y_p$ in $G$, and if $h_r, h_p$ are nonadjacent in $S^+$ then $\phi(h_r)$ is nonadjacent to every vertex
of $Y_p$ in $G$.
\end{itemize}
We call the list of sets $Y_{0}\LL Y_{p-1}$ a choice of {\em $t$-minions} for $\phi$, and $F$ the {\em image} of $\phi$.
We deduce:
\\
\\
(1) {\em No $1$-map is $2$-general.}
\\
\\
Suppose that $\phi$ is a $2$-general $1$-map, and let $Y_{0}$ be a choice of the (single) $t$-minion.
Since $\mu(Y_{0})>1$, $G[Y_{0}]$ is not $\mathcal{H}_2$-free; but there are no edges between $Y_{0}$ and $V(F)$,
and $F$ is isomorphic to $S$, contradicting that $G$ is $\mathcal{H}_1\uplus\mathcal{H}_2$-free. This proves (1).
\\
\\
(2) {\em For $1\le q<s$ and $t\ge 2$, if no $q$-map is $t$-general then no
$(q+1)$-map is $t'$-general, where 
$t'=2+(t-1)(k-2)$.
Consequently no $s$-map is $k^{s-1}$-general.
}
\\
\\
Suppose that $\phi$ is a $t'$-general $(q+1)$-map, let $Y_0\LL Y_{q}$ be some choice of $t$-minions for $\phi$, 
and let $F$ be its image. Let $v\in Y_q$. If we extend $\phi$ by mapping $h_q$ to $v$, we obtain a $q$-map, which therefore
is not $t$-general. Consequently, either:
\begin{itemize}
\item for some $p\in \{0\LL q-1\}$, $h_p,h_q$ are adjacent in $S^+$, and the set of vertices in
$Y_q$ adjacent to $v$ in $G$ has measure at most $t-1$; or
\item for some $p\in \{0\LL q-1\}$, $h_p,h_q$ are nonadjacent in $S^+$, and the set of vertices in
$Y_q$ nonadjacent to $v$ in $G$ has measure at most $t-1$.
\end{itemize}
In either case we say $p$ is a {\em problem} for $v$. For $0\le p<q$, let $V_p$ be the set of $v\in Y_q$ such that $p$
is a problem for $v$. Thus $V_0\cupcup V_{q-1}=Y_q$, and since 
$$\mu(Y_q)\ge t'\ge (t-1)|S|\ge |S|=s> q,$$ 
it follows that $\mu(V_p)>1$
for some $p\in \{0\LL q-1\}$. 

Suppose that $h_p,h_q$ are adjacent in $S^+$. Since $\mu(V_p)>1$, $G[V_p]$ contains a copy $F$ of some member $H_1\in \mathcal{H}_1$.
For each $v\in V(F)$, the set of vertices in $Y_q$ adjacent to $v$ has measure at most $t-1$,
and so the set of vertices in $Y_q$ with a neighbour in $V(F)$ has measure at most $|F|(t-1)=|H_1|(t-1)\le (k-2)(t-1)$.
Since $\mu(Y_q)\ge t'$, the set of vertices $X$ in $V(F_1)$ with no neighbour in $V(F_1)$ has measure at least $t'-(k-2)(t-1)= 2$
and so contains a copy of some member of $\mathcal{H}_2$, contradicting that $G$ is $\mathcal{H}_1\uplus\mathcal{H}_2$-free.

Now suppose that $h_p,h_q$ are nonadjacent in $S^+$. Since $\mu(V_p)>1$, $G[V_p]$ contains a copy $F$ of some member $J_1\in \mathcal{J}_1$.
For each $v\in V(F)$, the set of vertices in $Y_q$ nonadjacent to $v$ has measure at most $t-1$,
and so the set of vertices in $Y_q$ with a nonneighbour in $V(F_1)$ has measure at most $|F|(t-1)\le (k-2)(t-1)$.
Since $\mu(Y_q)\ge t'$, the set of vertices $X$ in $V(F_1)$ that are adjacent to every vertex in $V(F_1)$ has measure at least 
$t'-(k-2)(t-1)\ge 2$
and so contains a copy of some member of $\mathcal{J}_2$, contradicting that $G$ is $\mathcal{J}_1 * \mathcal{J}_2$-free. 
This proves the first assertion.

For the second assertion, (1) and the first assertion imply that no $2$-map is $k$-general. 
Also the first assertion implies that for $1\le q<s$ and $t\ge 1$, if no $q$-map is $t$-general then
no $(q+1)$-map is $kt$-general. From (1), it follows that for $1\le q\le s$ and $t\ge 1$, no $q$-map is $k^{q-1}$-general, and
in particular, no $s$-map is $k^{s-1}$-general.
This proves (2).

\bigskip

For each $v\in V(G)$, let $N(v)$ be the set of neighbours of $v$ in $G$, and let $M(v)$ be the set of its non-neighbours; so $M(v),N(v)$
have union $V(G)\setminus \{v\}$.
\\
\\
(3) {\em For every vertex $v\in V(G)$, either $\mu(N(v))< sk^{s-1}$ or $\mu(M(v))< sk^{s-1}$.}
\\
\\
Let $h_s$ have $a$ neighbours and $b$ non-neighbours in $S^+$. Thus $a\le s-1$ and $b\le s$.
Let $\phi$ be the function mapping $h_s$ to $v$. Thus $\phi$ is an $s$-map, and so is not $k^{s-1}$-general. Consequently,
either $N(v)$ cannot be partitioned into $a$ sets each with measure at least $k^{s-1}$, or $M(v)$ 
cannot be partitioned into $b$ sets each with measure at least $k^{s-1}$. 
But a minimal subset of $N(v)$ with measure at least $k^{s-1}$ has measure 
$k^{s-1}$ (since singleton sets are good); and so, by successively extracting minimal sets
with measure at least $k^{s-1}$, we deduce that either $\mu(N(v))< ak^{s-1}$ , or $\mu(M(v))< bk^{s-1}$.
In particular, one of $N(v), M(v)$ has measure less than $sk^{s-1}$.
This proves (3).

\bigskip

Let $A$ be the set of $v\in V(G)$ such that $N(v)$ has measure at most $sk^{s-1}$, and let $B$
be the set of $v\in V(G)$ such that $M(v)$ has measure at most $sk^{s-1}$.  If $\mu(A)>1$, then
$G[A]$ contains a copy $F$ of some member of $\mathcal{H}_1$, and the set of vertices of $G$ with a neighbour in $V(F)$ has measure at most
$sk^{s-1}|F|\le sk^s$; and since the set of vertices with no neighbours in $V(F)$ has measure at most 1 (because 
it induces an $\mathcal{H}_2$-free subgraph) and $V(F)$ has measure at most 2, it follows that $V(G)$ has measure at most
$sk^s+3$. If $\mu(B)>1$, then $G[B]$ contains a copy of some member of $\mathcal{J}_1$, and it follows similarly that 
$\mu(V(G))\le sk^s+3$.
If neither of these, then $\mu(V(G))\le \mu(A)+\mu(B)\le 2$. In each case the theorem holds. This proves \ref{pairsthm}.~\bbox

Let us deduce \ref{comps}, which is implied by the following:
\begin{thm}\label{union}
Let $H$ be the disjoint union of $H_1, H_2$, and suppose that $H_1,H_2$ are both multibounding. Then $H$ is multibounding.
\end{thm}
\Proof We need to show that for every $d\ge 1$, there 
is a polynomial $f_d(t)$ such that, for all $t\ge 1$,
if a graph $G$ is $H$-free and does not contain $K_d(t)$ as a subgraph, then $\chi(G)\le f_d(t)$. We show this by induction on $d$.
If $d=1$, the statement is trivial, taking $f_1(t)=t$, so we may assume that $d>1$, and the statement holds for all $d'<d$.
For $i = 1,2$, since $H_i$ is multibounding, there is a polynomial $g_i(t)$
such that for all $t\ge 1$,
if a graph $G$ is $H_i$-free and does not contain $K_d(t)$ as a subgraph, then $\chi(G)\le g_i(t)$. Let $g(t)$
be a polynomial such that $g_1(t), g_2(t)\le g(t)$, and $f_{d'}(t)\le g(t)$ for all $d'$ with $1\le d'<d$.

Let $k=2+|H_1|+d_1t$, and $s=|H_1|$,  and let $f_d(t)=(sk^s+3)g(t)$. We will show that for all $t\ge 1$,
if a graph $G$ is $H$-free and does not contain $K_d(t)$ as a subgraph, then $\chi(G)\le f_d(t)$. 

We may assume that at least one of $H_1,H_2$ has more than one vertex.
Let $\mathcal{H}_i=\{H_i\}$ for $i = 1,2$. Thus $G$ is $\mathcal{H}_1\uplus\mathcal{H}_2$-free.
Choose integers $d_1,d_2>0$ with $d_1+d_2=d$. For $i = 1,2$, let $\mathcal{J}_i$ be the class of all graphs with exactly
 $d_it$ vertices that have $K_{d_i}(t)$
as a subgraph. Since $G$ does not contain $K_d^t$ as a subgraph, it is $\mathcal{J}_1 * \mathcal{J}_2$-free,
and so we can apply \ref{pairsthm}. 
By \ref{pairsthm}, it follows that 
$V(G)$ can be partitioned into $n=sk^s+3$ sets $X_1\LL X_n$, such that for $1\le i\le n$, $G[X_i]$ is either
$\mathcal{H}_1$-free or $\mathcal{H}_2$-free or  $\mathcal{J}_1$-free or  $\mathcal{J}_2$-free. If $G[X_i]$ is $\mathcal{H}_1$-free,
then $\chi(G[X_i])\le g_1(t)\le g(t)$, and similarly if $G[X_i]$ is $\mathcal{H}_2$-free,
then $\chi(G[X_i])\le g(t)$. If $G[X_i]$ is $\mathcal{J}_1$-free, then it does not contain $K_{d_1}(t)$ as a subgraph,
and so $\chi(G[X_i])\le f_{d_1}(t)\le g(t)$, and similarly if 
$G[X_i]$ is $\mathcal{J}_2$-free, then $\chi(G[X_i])\le g(t)$. So in each case, $\chi(G[X_i])\le g(t)$, and hence 
$\chi(G)\le ng(t)=f_d(t)$. This proves \ref{comps}.~\bbox

While we are on the subject, let us give another nice application of \ref{pairsthm}. It is known that
if every component of a forest $H$ satisfies the Gy\'arf\'as-Sumner conjecture \ref{GSconj} then $H$ also satisfies it. 
It is not known whether, if every component of
a forest $H$ satisfies \ref{esperet}, then $H$ necessarily satisfies \ref{esperet} (see~\cite{poly6}); but this is true for a related property.
Let us say a graph $H$ is 
{\em near-Esperet} if there exists $c>0$ such that for all $d>0$, every $H$-free graph $G$ with $\omega(G)<d$
satisfies $\chi(G)\le d^{c \log d}$. (Logarithms are to base two throughout.)
For instance, the theorem of~\cite{poly5} says that the five-vertex path is near-Esperet. Here, we prove that if every component of $H$
is near-Esperet then $H$ is near-Esperet, because of the following:
\begin{thm}\label{nearesp}
Let $H$ be the disjoint union of $H_1,H_2$, and let $H_1,H_2$ be near-Esperet. Then $H$ is near-Esperet.
\end{thm}
\Proof
For $i = 1,2$
we know that $H_i$ is near-Esperet, so there exists $c_i>0$ such that for all $d>0$, every $\{H_i, K_d\}$-free graph $G$ satisfies
$\chi(G)\le d^{c_i\log d}$. Let $c_0=\max(c_1,c_2)$. Choose $b\ge 0$ such that $|H_1|(2+d+|H_1|)^{|H_1|}+3\le d^b$
for all $d\ge 2$.
Define $c=\max (b+c_0, b/\log (3/2))$.
We will
show by induction on $d$ that for all $d\ge 0$,  every $\{H,K_d\}$-free graph $G$ 
satisfies $\chi(G)\le d^{c \log d}$.
If $d\le 2$ then every $K_d$-free graph $G$ satisfies $\chi(G)\le 1$, and the claim holds since $c\ge 0$; so we may assume that 
$d\ge 2$, and the claim holds for all $d'<d$. (In fact we could start the induction at any convenient value of $d$, since we know
that $H_1\cup H_2$ satisfies \ref{GSconj}.)
We may assume that at least one of $H_1,H_2$ has more than one vertex. Choose integers $d_1,d_2\le 2d/3$ with $d_1+d_2=d$, and let $J_i$ be the complete graph with $d_i$ vertices for $i = 1,2$.

Now let $G$ be $H$-free, with clique number less than $d$. Then $G$ does not contain the complete join of $J_1,J_2$, so by 
\ref{pairsthm}, taking $\mathcal{H}_i=\{H_i\}$ and $\mathcal{J}_i=\{J_i\}$ for $i = 1,2$, there is a partition of $V(G)$
into $n$ sets $X_1\LL X_n$, where $k=2+\max(d_1,|H_1|)$ and $n= |H_1|k^{|H_1|}+3$, such that for $1\le i\le n$, $G[X_i]$ is either
$H_1$-free or $H_2$-free or $J_1$-free or $J_2$-free. If $G[X_i]$ is $H_1$-free then $\chi(X_i)\le d^{c_1\log d}\le d^{c_0\log d}$, and 
similarly if $G[X_i]$ is $H_2$-free then $\chi(X_i)\le d^{c_0\log d}$. If $G[X_i]$ is $J_1$-free then from the inductive hypothesis,
$$\chi(X_i)\le d_1^{c\log d_1}\le (2d/3)^{c \log(2d/3)},$$ 
and similarly if  $G[X_i]$ is $J_2$-free then $\chi(X_i)\le (2d/3)^{c \log(2d/3)}$. Since $n\le d^b$,
it follows that 
$$\chi(G)\le d^b\max\left(d^{c_0\log d},(2d/3)^{c \log(2d/3)}\right).$$
Thus, to show that $\chi(G)\le d^{c\log d}$,  it suffices to check that 
$$d^{c\log d}\ge d^b d^{c_0\log d}$$
and 
$$d^{c\log d}\ge d^b (2d/3)^{c \log(2d/3)}.$$
The first says (taking logarithms)
$$c(\log d)^2\ge b\log d + c_0(\log d)^2,$$
and this is true since $b\log d\le b(\log d)^2$ and $c\ge b+c_0$.
The second says (taking logarithms)
$$c (\log d)^2 \ge b \log d+ c (\log (2d/3))^2.$$
Let us write $\log d = x$ and $\log (3/2) = \vare$; then we need to show that
$$cx^2\ge bx+ c (x-\vare)^2,$$
that is, $2c\vare x \ge bx+c\vare^2$. But this is true since
$c\vare x\ge bx$ and $c\vare x\ge c\vare^2$ (because $d>3/2$).
This proves \ref{nearesp}.~\bbox

It is proved in~\cite{poly4} that the five-vertex path $P_5$ is near-Esperet; so an application of \ref{nearesp} yields 
that every 
forest that is the disjoint union of copies of $P_5$ is also near-Esperet.

\end{document}